\documentclass[12pt,letterpaper]{amsart}
\usepackage{euler, epic,eepic,latexsym, amssymb, amscd, amsfonts, xypic, floatflt}
\input xy
\xyoption{all}

 \newlength{\baseunit}               
 \newcount{\numlines}                
\newcommand{\tpoint}[1]{\vspace{3mm}\par \noindent \refstepcounter{subsection}{\bf \thesubsection.} 
  {\em #1. ---} }
\newcommand{\epoint}[1]{\vspace{3mm}\par \noindent \refstepcounter{subsection}{\bf \thesubsection.} 
  {\em #1.} }
\newcommand{\bpoint}[1]{\vspace{3mm}\par \noindent \refstepcounter{subsection}{\bf \thesubsection.} 
  {\bf #1.} }

\newcommand{\bpf}{\noindent {\em Proof.  }}
\newcommand{\epf}{\qed \vspace{+10pt}}

\newcommand{\zed}{\mathbb{Z}}
\newcommand{\Z}{\mathbb{Z}}
\newcommand{\A}{\mathbb{A}}

\newcommand{\C}{\mathbb{C}}
\newcommand{\F}{\mathbb{F}}

\newcommand{\M}{M}

\newcommand{\proj}{\mathbb P}
\newcommand{\tS}{\tilde{S}}
\newcommand{\oh}{{\mathcal{O}}}
\newcommand{\cI}{{\mathcal{I}}}
\newcommand{\cJ}{{\mathcal{J}}}
\newcommand{\cK}{{\mathcal{K}}}

\renewcommand{\cL}{{\mathcal{L}}}
\newcommand{\cm}{{\mathcal{M}}}
\newcommand{\cu}{{\mathcal{U}}}
\newcommand{\cX}{{\mathcal{X}}}
\newcommand{\cY}{{\mathcal{Y}}}
\newcommand{\al}{\alpha}

\newcommand{\be}{\beta}
\newcommand{\ga}{\gamma}

\newcommand{\si}{\sigma}

\newcommand{\Hom}{\operatorname{Hom}}
\newcommand{\Ext}{\operatorname{Ext}}
\newcommand{\sHom}{\operatorname{\mathcal{H}om}}
\newcommand{\sExt}{\operatorname{\mathcal{E}xt}}
\newcommand{\Hilb}{\operatorname{Hilb}}
\newcommand{\Pic}{\operatorname{Pic}}
\newcommand{\Def}{\operatorname{Def}}
\newcommand{\Spec}{\operatorname{Spec}}

\newcommand{\Aut}{\operatorname{Aut}}
\newcommand{\codim}{\operatorname{codim}}
\newcommand{\Div}{\operatorname{Div}}

\newcommand{\cited}{}

\newcommand{\remind}[1]{{\bf[#1]}}
\newcommand{\notation}[1]{}
\renewcommand{\remind}[1]{{}}
\newcommand{\lremind}[1]{{}}

\newcommand{\secretnote}[1]{}

\begin{document}
\pagestyle{plain}
\title{{\large {Murphy's Law in algebraic geometry:  \\
 Badly-behaved deformation spaces}}
}
\author{Ravi Vakil}
\address{Dept. of Mathematics, Stanford University, Stanford CA~94305--2125}
\email{vakil@math.stanford.edu}
\thanks{Partially supported by NSF CAREER/PECASE Grant DMS--0228011, and 
an Alfred P. Sloan Research Fellowship.
%
}
\date{Sunday, November 21, 2004.}
\begin{abstract}
We consider the question: ``How bad can the deformation space of an
object be?'' The answer seems to be: ``Unless there is some a priori
reason otherwise, the deformation space may be as bad as
possible.''  We show this for a number of important moduli spaces.

More precisely, every singularity of finite type over $\zed$ (up
to smooth parameters) appears on: the Hilbert scheme of curves in
projective space; and the moduli spaces of smooth projective general-type
surfaces (or higher-dimensional varieties), plane curves with nodes
and cusps, stable sheaves, isolated threefold singularities, and more.
The objects themselves are not pathological, and are in fact as nice
as can be: the curves are smooth, the surfaces have very ample
canonical bundle, the stable sheaves are torsion-free of rank $1$, the
singularities are normal and Cohen-Macaulay, etc.  This justifies
Mumford's philosophy that even moduli spaces of well-behaved objects
should be arbitrarily bad unless there is an a priori reason
otherwise.

Thus one can construct a smooth curve in projective space whose
deformation space has any given number of components, each with any
given singularity type, with any given non-reduced behavior along
various associated subschemes.  Similarly one can give a surface over
$\F_p$ that lifts to $\Z/p^7$ but not $\Z/p^8$.  (Of course the results hold
in the holomorphic category as well.)

It is usually difficult to compute deformation spaces directly from
obstruction theories.  We circumvent this by relating them to
more tractable deformation spaces
via smooth
morphisms.  The essential starting point is
Mn\"ev's Universality Theorem.
\end{abstract}
\maketitle
\tableofcontents

{\parskip=12pt 

\begin{center}
\begin{minipage}{4in}
{\em 
\begin{verse}
The best-laid schemes o' mice an' men  \\
Gang aft agley \\
An' lea'e us nought but grief an' pain \\
For promis'd joy! \\ 
\quad
\quad
\quad
--- Robert Burns, ``To a Mouse''
\end{verse}
}
\end{minipage} 
\end{center}

\section{Introduction}
Define an equivalence relation on pointed schemes
generated by:  If $(X,p) \rightarrow (Y,q)$ is a smooth morphism,
then $(X,p) \sim (Y,q)$.
We call the equivalence classes {\em singularity types},
and will call pointed schemes {\em singularities} (even if the
point is regular). 
We say that {\em Murphy's Law}
holds for a moduli space if every singularity type of finite type
over $\mathbb{Z}$  appears on that moduli space.
Although our methods are algebraic, our arguments all work in the
holomorphic category. 


\tpoint{Main Theorem}  \label{mainthm}
{\em The following moduli spaces satisfy Murphy's Law. \lremind{mainthm}
\begin{enumerate}
\item[{\bf \M 1a.}]  the Hilbert scheme of nonsingular curves in projective
space 
\item[{\bf \M 1b.}]  the moduli space of maps of smooth curves to
projective space (and hence Kontsevich's moduli space of maps)
\item[{\bf \M 1c.}] 
$\mathcal{G}^r_d$ \cite[p.~5]{hm}, the space of  curves
with the data of a linear system of degree $d$ and projective dimension $r$
\item[{\bf \M 2a.}]  the versal deformation spaces of smooth surfaces
  (with very ample canonical bundle) 
\item[{\bf \M 2b.}] the fine moduli stack of smooth
surfaces with very ample canonical bundle and reduced automorphism group
\item[{\bf \M 2c.}] the coarse moduli space of smooth surfaces
with very ample canonical bundle
\item[{\bf \M 2d.}]  the Hilbert scheme of nonsingular surfaces in 
$\mathbb{P}^5$, and the Hilbert scheme of surfaces in $\mathbb{P}^4$
\item[{\bf \M 3a--c.}]  more generally, the versal deformation spaces 
and fine and coarse moduli spaces of
  smooth $n$-folds ($n>1$) with very ample canonical bundle and 
reduced automorphism group (as in {\bf 2a--c})
\item[{\bf \M 4.}]  the Chow variety of nonsingular curves in projective
space, and of nonsingular surfaces in $\proj^5$, 
allowing only seminormal singularities in the definition of Murphy's
Law
(recall that the Chow variety is seminormal \cite[Theorem~3.21]{kollar})
\item[{\bf \M 5a.}]  branched covers of $\mathbb{P}^2$ with only simple branching
(nodes and cusps), in characteristic not $2$ or $3$
\item[{\bf \M 5b.}]  the ``Severi variety'' of plane curves with 
a fixed numbers of nodes and cusps, in characteristic not $2$ or $3$
\item[{\bf \M 6.}] the moduli space of stable sheaves  \cite{simpson}
\item[{\bf \M 7.}] the versal 
deformation spaces of isolated  normal Cohen-Macaulay
threefold singularities
\end{enumerate}}

The meaning of Murphy's Law for versal deformation spaces is the
obvious one.  We should say a few words on why certain moduli spaces
exist. {\bf 1b:} Although one usually discusses Kontsevich's moduli
space of stable maps in characteristic $0$, one may as well define the
moduli space of maps from nodal curves to projective space, with
reduced automorphism group, over $\Spec \Z$; this is a Deligne-Mumford
stack, by a standard generalization of the construction of \cite{fp}.
(It is not proper!)  {\bf 2b} and {\bf 3b}: \cite[p.~182--3]{artin}
shows existence for surfaces, and the argument applies verbatim in
higher dimension.  The stack is Deligne-Mumford, locally of finite
type.  {\bf 2c} and {\bf 3c}: \cite[Theorem~1.8]{kollarq} shows that
there is an algebraic space coarsely representing these moduli
functors.  (For surfaces {\bf 2c}, there is even a coarse moduli
(algebraic) space of canonical models of surfaces of general type
\cite[Theorem~1.7]{kollarq}.)

To obtain results over other bases (such as algebraically closed
fields such as $\C$), note that most moduli spaces above behave well with
respect to base change; hence any singularity obtained by base change
from a finite type singularity over $\Z$ may appear.  
Clearly, no other singularity may appear.
Indeed, any moduli (pseudo-)functor admitting a smooth cover by a scheme
locally of finite type over $\Spec \Z$ necessarily only has
singularities of this sort.
(For example, the singularity   \lremind{piinthesky}
\begin{equation}
\label{piinthesky}
xy(y-x)(y- \pi x) =0
\end{equation} 
in $\C^2$ may not appear as such a deformation space.)  This leads to some
natural questions, such as: does there exist an isolated complex
singularity whose deformation space is equivalent to \eqref{piinthesky}?  What
if the singularity is required to be algebraic?  Does there exist a
compact complex manifold whose deformation space has such singularity
type?  What if the manifold were required to be projective? 
Before this project, we would have believed the answer could be ``yes'',
but paradoxically Murphy's Law leads us now to expect that the answer is
``no''.  It would be very interesting to have any example of a
non-pathological object (e.g.\ isolated complex algebraic singularity,
complex projective manifold, or even non-algebraic examples) with deformation space not equivalent to
one of finite type over $\Spec \Z$.

\remind{Old note:  Brian's inserts are in his e-mail, and in Mnov1404.tex.}


\bpoint{Philosophy}
To be explicit about why these results may be surprising: one can
construct a smooth curve in projective space whose deformation space
has any given number of components, each with any given singularity
type, with any given non-reduced behavior along various associated
subschemes.  Similarly, one can give a smooth surface of general type
in characteristic $p$ that lifts to $\Z/p^7$ but not to $\Z/p^8$.


We next give some philosophical comments, which motivated this result.
The history sketched in Section~\ref{history} also provided motivation.

The moral of Theorem~\ref{mainthm} is as follows.  We know that some
moduli spaces of interest are ``well-behaved'' (e.g.\ equidimensional,
having at worst finite quotient singularities, etc.), often because
they are constructed as Geometric Invariant Theory quotients of smooth
spaces, e.g.\ the moduli space of curves, the moduli space of vector
bundles on a curve, the moduli space of branched covers of
$\mathbb{P}^1$ (the Hurwitz scheme, or space of admissible or twisted
covers), the Picard variety, the Hilbert scheme of divisors on
projective space, the Severi variety of plane curves with a prescribed
number of nodes, the moduli space of abelian varieties (notably \cite{no}),
etc.  In other cases, there has been some effort to try to
bound how ``bad'' the singularities can get.  Theorem~\ref{mainthm} in
essence states that these spaces can be arbitrarily singular, and
gives a means of constructing an example where any given behavior
happens.

Murphy's law suggests that unless there is some natural reason for the
space to be well-behaved, it will be arbitrarily badly behaved.  For
example, arithmetically Cohen-Macaulay surfaces in $\proj^4$ are
always unobstructed \cite{e}; but surfaces in general in $\proj^4$ can
have arbitrarily bad deformations (by {\bf \M 2d}).  Other
examples are given in the following table.

\begin{center}
\begin{tabular}{|l|l|} \hline
{\em Well-behaved moduli space} & {\em Badly-behaved moduli space} \\ \hline 
curves & surfaces  (by {\bf \M 2b--c}) \\
branched covers of $\proj^1$ (e.g.\ \cite[Theorem~1.53]{hm}) &
  branched covers of $\proj^2$ (by {\bf \M 5a}) \\
surfaces in $\proj^3$ & surfaces in $\proj^4$ (by {\bf \M 2d}) \\
Picard variety over the moduli space of curves &
its subscheme $\mathcal{G}^r_d$  (by {\bf \M 1c}) \\
Severi variety of nodal   plane curves  &
Severi variety of nodal and \\
\quad  (e.g.\ \cite[Theorem~1.49]{hm}) &
\quad cuspidal plane curves (by {\bf \M 5b}) \\ \hline
\end{tabular}
\end{center}


Furthermore, our experience and intuition tells us that pathologies of
moduli spaces occur on the boundary, and that moduli spaces of
``good'' objects are also ``good''.  Murphy's Law shows that this
intuition is incorrect; we should expect pathologies even where the
objects being parameterized seem harmless.  Kodaira says ``The theory
of deformation was at first an experimental science''
\cite[p.~259]{kodairabk}.  This result shows that our intuition is
flawed because it is based on experimental knowledge of a very small
part of the moduli spaces we are interested in; it supports Mumford's
philosophy that pathologies are the rule rather than the exception.
Alternatively, from the point of view of A.~Vershik, this result
states that the ``universality'' philosophy
(e.g.~\cite[Section~7]{vershik}) applies widely in algebraic geometry.


As a side comment, Theorem~\ref{mainthm} indicates that one cannot
hope to desingularize the moduli space of surfaces, or any other moduli
space satisfying Murphy's Law, by adding additional structure; this
would imply a resolution of singularities.  (Hence the program for
desingularization of the space of stable maps informally proposed by
some authors cannot succeed.  However, see \cite{vz} for success in
genus $1$.)

\epoint{Notation} Let $\Def$ denote the versal or Kuranishi
deformation space (not the space of first-order deformations).  The
object being deformed will be clear from the context.

\bpoint{Acknowledgments} I am indebted to A. J. de Jong and S. Billey
for discussions that led to these ideas.  I am grateful to the
organizers and participants in the 2004 Oberwolfach workshop on
Classical Algebraic Geometry for many comments.  I thank B. Shapiro in
particular for pointing out that Theorem \ref{mut} was first proved by
Mn\"{e}v.  I thank F.~Catanese, R.~Thomas, J.~Wahl, M.~van Opstall,
R.~Pardini, M.~Manetti, B.~Conrad, B.~Hassett and S.~Kov\'{a}cs for sharing their expertise.
Significant improvements to this paper are due to them.  I also thank
W.~Fulton and A.~Vershik.

\section{History, and further questions}
\label{history}

\bpoint{Hilbert schemes} \lremind{history}The 
motivation for both the equivalence relation $\sim$ and the
terminology ``Murphy's Law'' comes from the folklore conjecture that
the Hilbert scheme ``satisfies Murphy's Law''.  

\tpoint{Law \cite[p.~18]{hm}}{\em 
There is no geometric possibility so horrible that it cannot
be found generically on some component of some Hilbert scheme.}

I am not sure of the origin of this philosophy, but it seems
reasonable to ascribe it to Mumford.  This traditional statement of
Murphy's Law is admittedly informal and imprecise (see the MathReview
\cite{hmmr}).  Clearly not every singularity appears
on the Hilbert scheme of projective space.  For example, the only
zero-dimensional Hilbert schemes are reduced points.  Allowing
``smooth equivalence classes'' of singularities seems the mildest way
of rescuing the law.

In his famous paper \cite{mumford}, Mumford described a component of the
Hilbert scheme of space curves that is everywhere nonreduced.  
Other examples of nonreduced components of the Hilbert scheme have
since been given \cite{gp, kleppe, ellia, mdp}.  Other pathologies relating
to the number of components of the Hilbert scheme of {\em smooth}
space curves were given by Ellia, Hirschowitz, and Mezzetti
\cite{ehm}, and by Fantechi and Pardini \cite{fp1}.  (The results of
the latter will be essential to our argument.)

Raynaud's example (see Section~\ref{multiple}) gives a component of a
Hilbert scheme of smooth surfaces which exists in characteristic $p$,
but does not lift to characteristic $0$ (by the standard methods of
Section~\ref{pf2d}).  Mohan Kumar, Peterson, and Rao \cite{mpr} give a
component of the Hilbert scheme of smooth surfaces in $\proj^4$ which
exists in characteristic $2$ but does not lift.  See \cite[Section~3]{elha}
for more on problems of lifting curves out of characteristic $p$.

Although the Hilbert scheme of projective spaces was suspected
to behave badly, other moduli spaces were believed (or hoped) to be
better-behaved.  We now discuss these.

\bpoint{Surfaces and higher-dimensional varieties}
\label{multiple}
\lremind{multiple}(See \cite{cataneseMAS} for an excellent overview of
the subject.)  The first example of an obstructed smooth variety was
due to Mumford, obtained by blowing up his curve in $\proj^3$
\cite[p.~643-4]{mumford}.  The first example of an obstructed surface
is due to Kas \cite{kas}.  Other examples were later given by Burns
and Wahl \cite{bw}, and later many others.  Horikawa \cite{horikawa},
Miranda \cite{miranda}, and Catanese \cite{cataneseNR} gave
examples of generically nonreduced components of the moduli space of
surfaces; in each case the surfaces did not have ample canonical
bundle, and this appeared to be a common explanation of this pathology
\cite[p.~294]{cataneseNR}.  Catanese conjectured that if $S$ is a
surface of general type with $q=0$ and $K_S$ ample, then the moduli
space $\mathfrak{M}(S)$ is smooth on an open dense set (\cite[p.~34,
69]{cataneseMAS}, \cite[p.~294]{cataneseNR}).  Theorem~\ref{mainthm}
{\bf \M 2b--c} gives a counterexample to this conjecture, and as Catanese
pointed out, even to the stronger conjecture where $K_S$ is very
ample.  Manetti gave an earlier counterexample in his thesis
\cite[Corollary~3.4]{manettithesis}; the added advantage of
{\bf \M 2} is that every (finite type) nonreduced structure
is shown to occur.

Catanese showed that the
moduli space of complex surfaces in a given homeomorphism class can
have arbitrarily many components of different dimension
\cite[Theorem~A]{catanese}, and asked if this were still true for
those in a given diffeomorphism class \cite[p.~485]{catanese};
Theorem~\ref{mainthm} {\bf \M 2b--c} answers this in the affirmative.
A prior answer was recently given by Catanese and B.~Wajnryb
\cite{cw}.  The added benefit of {\bf \M 2} is that all possibilities
are shown to occur.

Serre gave the first example of a projective variety that could not
be lifted to characteristic $0$ \cite{serre}. Raynaud
gave the first example of such a surface \cite{raynaud};
W. Lang gave more \cite{lang}.

\bpoint{Plane curves with nodes and cusps} If $C$ is a reduced complex
plane curve, the classical question of ``completeness of the
characteristic linear series'' asks (in modern language) if an
appropriate equisingular moduli space is smooth.  Severi proved this
is true if $C$ has only nodes (\cite{severi}, see also
\cite[Section~VIII.4]{zariski}), and  asserted this 
if $C$ has nodes and cusps \cite{wmr}.
(See \cite[p.~116-7 and Section~VIII]{zariski} for motivation for the study of
nodal and cuspidal plane curves.)  
It was later realized that Severi's
assertion was unjustified.
Enriques tried repeatedly 
to show that such curves were unobstructed \cite[p.~51]{cataneseMAS};
Zariski also raised this question 
\cite[p.~221]{zariski}.  The first counterexample was given by Wahl 
\cite[Section~3.6]{wahl}, and another was given by Luengo
\cite{luengo}.  Theorem~\ref{mainthm} {\bf \M 5b} shows that Severi
was in some sense ``maximally wrong''.

\bpoint{Stable coherent sheaves} The moduli space of stable coherent
sheaves is due to Simpson \cite{simpson}.  Our example is in fact a
torsion-free sheaf on $\proj^5$; the theory of the moduli of
torsion-free sheaves was developed earlier by Maruyama
\cite{maruyama}, building on Gieseker's work in the surface case
\cite{gieseker}.

\bpoint{Singularities} The theory of deformations of singularities is
too large to summarize  here.  We point out however that it was
already established by Burns and Wahl \cite{bw} that such deformation
spaces can be bad, although not this pathological.

\bpoint{Further questions} Theorem~\ref{mainthm}, and the philosophy
and history given above, beg further questions.  Do deformations of
surface singularities (say isolated and Cohen-Macaulay) satisfy Murphy's
Law?  How about the Hilbert scheme of curves in $\proj^3$?  
The Hilbert scheme of points on a smooth threefold?
The moduli
of vector bundles on smooth surfaces?  
Can the extra dimensions
allowed in the definition of {\em type} be excised, i.e.\ can
``smooth'' be replaced by ``\'{e}tale'' in the definition of type?
(As observed above, this is not possible for the Hilbert scheme.)
Catanese asks if Murphy's law for surfaces is still true if we require
not only that the surface has very ample canonical bundle, but also
that the canonical embedding is cut out by quadrics.  {\em Conjecture:} for any given $p$,
the surfaces whose canonical divisor induces an embedding satisfying
property $N_p$ satisfy Murphy's Law.  One might hope
that the constructions given in Section~\ref{abeliancovers} suffice by
taking the divisor class $A$ (Sections~\ref{p2} and~\ref{p3}) to be
sufficiently ample.  The case $p=1$ would give an affirmative answer
to Catanese's question.

\section{The starting point:
Mn\"{e}v's Universality Theorem}

\label{mnev}
\remind{Old note:  
Can I cut the following paragraph?:  [Given two moduli spaces
$\cm_1$ and $\cm_2$, and points $p_1$, $p_2$ of each we say the
deformations of $p_1$ are the same (respectively of the same type) as
those of $p_2$ if there are Zariski-open neighborhoods $\cu_1 \subset
\cm_1$ and $\cu_2 \subset \cm_2$ and an isomorphism (respectively
smooth morphism) $\cu_1 \rightarrow \cu_2$.]}

We will prove
Theorem~\ref{mainthm} by drawing connections among various moduli
spaces, taking as a starting point a remarkable result of Mn\"{e}v.
Define an {\em incidence scheme of points and lines in
  $\mathbb{P}^2$,} a locally closed subscheme of $(\mathbb{P}^2)^m
\times (\mathbb{P}^{2*})^n = \{ p_1, \dots, p_m, l_1, \dots, l_n \}$
parameterizing $m \geq 4$ marked points and $n$ marked lines, as follows.

\noindent $\bullet$ $p_1 = [1;0;0]$, $p_2 = [0;1;0]$, $p_3 = [0;0;1]$,
$p_4 = [1;1;1]$.
\newline $\bullet$ We are given
some specified incidences: For each pair $(p_i, l_j)$, either $p_i$ is
required to lie on $l_j$, or $p_i$ is required not to lie on $l_j$.
\newline $\bullet$ The marked points are required to be distinct, and the
marked lines are required to be distinct.  
\newline $\bullet$ Given any two
marked lines, there a marked point required to be on both of them.  
\newline
$\bullet$ Each marked line contains at least three marked points.

\tpoint{Theorem (Mn\"{e}v)} {\em Every singularity type appears on
  some incidence scheme.} \label{mut}

This is a special case of Mn\"{e}v's Universality Theorem \cite{mnev1,
  mnev2}.  A short proof is given by Lafforgue in
\cite[Th\'{e}or\`{e}me~1.14]{lafforgue}.  Lafforgue's construction
does not necessarily satisfy the first, fourth and fifth requirements
of an incidence scheme, but they can be satisfied by adding more
points.  (The only subtlety in adding these extra points is verifying
that in the configuration constructed by Lafforgue, no three lines
pass through the same point unless required to by the construction.)
{\em Caution:} Other expositions of Mn\"{e}v's theorem do not prove
the result scheme-theoretically, only ``variety-theoretically,'' as
this is all that is needed for most purposes. 

For the rest of the paper fix a singularity type.  By Mn\"{e}v's
Theorem~\ref{mut}, there is an incidence scheme exhibiting this
singularity type at a certain configuration $\{ p_1, \dots, p_m,
l_1, \dots, l_n \}$.  
Consider the surface $S$ that is the blow-up of $\mathbb{P}^2$ at the
points $p_i$.  Let $C$ be the proper transform of the union of the
$l_j$, so $C$ is a smooth curve (a union of $\proj^1$'s).  
This induces a morphism from the incidence scheme to the 
moduli space of surfaces with marked smooth divisors.

\tpoint{Proposition} {\em This morphism is \'{e}tale at $(\proj^2, \{
  p_i \}, \{l_j \}) \mapsto (S,C)$.}

Thus the singularity at $(\proj^2, \{ p_i \}, \{l_j \})$ has the
same type as the moduli space of surfaces with marked smooth divisor
at $(S,C)$.

\bpf
We will produce an \'{e}tale-local inverse near $(S,C)$.
Consider a deformation of $(S,C)$:\lremind{murder}
\begin{equation}
\label{murder}
\xymatrix{
(S,C) \ar[d] \ar @{^{(}->}[r] & (\mathcal{S}, \mathcal{C}) \ar[d] \\
\text{pt} \ar  @{^{(}->}[r] & B.
}
\end{equation}
Pull back to an \'{e}tale neighborhood of $\text{pt}$ so that the
components of $C$ are labeled.  The Hilbert scheme of $(-1)$-curves is
\'{e}tale over the base.  (I am not aware of the first reference
for this well-known  fact.  It follows for example
from the exact sequence of \cite{ran} --- see e.g.\ the proof
of \cite[Theorem~3.2]{ran} --- which specializes to give a natural bijection
between the deformations, respectively obstructions,
of $S$ and $\proj^1 \rightarrow S$.  The proof in the holomorphic
category is due to Kodaira \cite[Theorem~3]{kodaira}.)

Let $E_i$ be the $(-1)$-curve corresponding to $p_i$.  Pull back to an
\'{e}tale neighborhood so that the points of the Hilbert scheme
corresponding to $E_i$ extend to sections (so there are divisors
$\mathcal{E}_i$ on the total space of the family that are
$(-1)$-curves on the fibers).  By abuse of notation, we use the same
notation \eqref{murder} for the resulting family.  By Castelnuovo's
criterion, $\mathcal{S}$ can be blown down along the $\mathcal{E}_i$
so that the resulting surface is smooth, with marked sections
extending $\{ p_i \}$.  (Again, ``Castelnuovo's criterion over
an Artin local scheme'' 
is presumably well-known to experts, but I am
unaware of a reference.  It follows by applying the ``usual''
Castelnuovo criterion over the closed point, and then using 
Theorem~\ref{wicked} to show that the blow-down ``deforms''.  
This is just the
old idea proved by  Horikawa 
in the smooth holomorphic category \cite{horikawa}.
Alternatively,
the proof of the usual
Castelnuovo's criterion, for example \cite[Theorem~V.5.7]{hartshorne},
can be extended.)

The central fiber is then $\proj^2$, so (as $\proj^2$ is rigid)
the family is locally trivial.  The marked points $p_1$, \dots, $p_4$ give a
canonical isomorphism with $\proj^2$.  (We may need to restrict to a smaller
neighborhood to ensure that these points are in general position.)  As
the components $\{ C_j \}$ of $C$ necessarily meet various $E_i$,
their images $\{ l_j \}$ necessarily pass through the necessary $p_i$.
\epf

\section{Abelian covers:  proof of
{\bf \M 2}}
\label{abeliancovers}
We use this intermediate moduli space of surfaces with marked divisors
to prove {\bf \M 2}, by connecting such marked surfaces to abelian
covers.  We use the theory of abelian covers developed by Catanese,
Pardini, Fantechi, and Manetti \cite{catanese, pardini, fp1, manetti}.
(Bidouble covers were introduced by Catanese.  Pardini developed the
general theory of abelian covers.  Key deformation-theoretic results
were established by Fantechi-Pardini and Manetti.)  Let $G =
(\zed/p)^3$, where $p=2$ or $3$ is prime to the characteristic of the
residue field of the singularity.  Let $G^\vee$ be the dual group, or
equivalently the group of characters.  Let $\langle \cdot, \cdot
\rangle: G \times G^\vee \rightarrow \zed/p$ be the pairing (after
choice of root of unity $\zeta$), which we extend to $\langle \cdot, \cdot
\rangle: G \times G^\vee \rightarrow \zed$ by requiring $\langle \si,
\chi \rangle \in \{ 0, \dots, p-1 \}$.  Suppose we have two maps $D: G
\rightarrow \Div(S)$, $L: G^\vee \rightarrow \Pic(S)$.  We say $(D,L)$
satisfies the {\em cover condition} \cite[Proposition~2.1]{pardini} if
$(D,L)$ satisfies $D_0=0$ and
$$
p  L_{\chi} =  \sum_\si  \langle \si,  \chi \rangle  D_{\si}
$$
for all $\si$, $\chi$.  (Equality is taken in $\Pic(S)$.)

\tpoint{Proposition (Pardini)} \label{pardinirocks} {\em Suppose $(D,L)$
  satisfies the cover condition, and suppose the $D_\si$ are
  nonsingular curves, no three meeting in a point, such that if
  $D_\si$ and $D_{\si'}$ meet then they are  transverse and $\si$
  and $\si'$ are linearly independent in $G$.  Then: \lremind{pardinirocks}
\begin{enumerate}
\item[(i)] There is a corresponding $G$-cover $\pi: \tS \rightarrow S$ with
branch divisor $D = \cup D_{\si}$. 
\item[(ii)] $\tS$ is nonsingular.
\item[(iii)] $\pi_* \oh_{\tS} = \oplus_{\chi} \oh_S(-L_{\chi})$.
\item[(iv)] $\pi_* \cK_{\tS} \cong \oplus_\chi \cK_S(L_\chi)$.
The Galois group $G$ acts on the left side in the obvious way;
it acts on the $\chi$-summand on the right by the character
$\chi$.
\end{enumerate}
}

Note for future reference that the branch divisor $D_{\si}$
corresponds to the subgroup of $G$ generated by $\si$.
(Note also that (iii) and (iv) are consistent with 
Serre duality on $\tS$.)

\bpf (i) is \cite[Proposition~2.1]{pardini}, (ii) is
\cite[Proposition~3.1]{pardini}, and (iii) is a consequence of Pardini's
construction \cite[(1.1)]{pardini}.
Pardini points out that (iv) is a special case of
duality for finite flat morphisms, see \cite{hartshorne} Exercises III.6.10
and Ex.~III.7.2.
(It also follows by a straightforward local calculation.
See \cite[p.~495]{catanese} for the analogous proof
for bidouble covers.  The generalization to abelian covers
is analogous to Pardini's proof of (iii).)
\epf

The next two examples apply to $(S,C)$ produced at the end of
Section~\ref{mnev}.  If the character of the residue field is $2$
(respectively $3$), then only Example~\ref{p3} (respectively \ref{p2})
applies; otherwise both apply.

\epoint{Key example: $p=2$} \label{p2} Fix $\si_0 \neq 0$ in $G$.  Let $A$ be a
sufficiently ample bundle such that $A \equiv C \pmod 2$.  Let
$D_{\si_0} = C$, $D_{0}=0$, and let $D_{\si}$ be a general section of
$A$ otherwise, such that $D_{\si'}$ meets $D_{\si''}$ transversely for
all $\si' \neq \si''$.  Let $L_0=0$, $L_\chi = 2A$ if $\langle \si_0, \chi
 \rangle = 0$ and $\chi \neq 0$, and $L_{\chi} = (3A+C)/2$ if
$\langle \si_0, \chi \rangle = 1$.  
(As $\Pic S$ is torsion-free, there is no ambiguity in the
phrase $(3A+C)/2$.)
It is straightforward
to verify that $(D,L)$ satisfies the hypotheses of Proposition~\ref{pardinirocks}.
\lremind{p2}

\epoint{Key example: $p=3$} \label{p3} Fix $\si_0 \neq 0$ in $G$, and
$\chi_0 \in G^\vee$ such that $\langle \si_0, \chi_0 \rangle = 1$.
Let $A$ be a sufficiently ample bundle such that $A \equiv C \pmod 3$.
Let $D_{\si_0} = C$, $D_{\si}$ be a general section of $A$ if $\langle
\si, \chi_0 \rangle = 1$ and $\si \neq \si_0$, and $D_{\si}=0$
otherwise.  Let
\begin{itemize}
\item $L_\chi = (8A+C)/3$ if $\langle \si_0, \chi \rangle = 1$
\item $L_0 = 0$
\item $L_\chi = 3A$ if $\langle \si_0, \chi \rangle = 0$ and $\chi \neq 0$
\item $L_{- \chi_0} = (16A+2C)/3$
\item $L_\chi = (7A+2C)/3$ if $\langle \si_0, \chi \rangle = 2$ 
and $\chi \neq - \chi_0$
\end{itemize}
It is straightforward
to verify that $(D,L)$ satisfies the hypotheses of Proposition~\ref{pardinirocks}
(note that if $\si \neq 0$, then at most one of $\{ D_{\si}, D_{-\si} \}$ 
is nonzero).


\tpoint{Theorem}
{\em  In Examples~\ref{p2} and~\ref{p3}, if $A$ is  sufficiently ample,
then $K_{\tS}$  is very ample.  In particular, $\tS$ is of general type,
and is its own canonical model.} \label{mufti} \lremind{mufti}

It is not hard to show that $K_{\tS}$ is {\em ample}:
$$
2 K_{\tS} = \pi^* \left( 2 K_S + \sum D_{\si} \right)
= \pi^* ( 2 K_S + C+ qA)$$
where $q=6$ if $p=2$ and $q=8$ if $p=3$.
If $A$ is  sufficiently ample, then
$2 K_S + \sum D_\si$ is ample, hence (as $\pi$ is finite)
$K_{\tS}$ is ample.
I am grateful to F.~Catanese for pointing out that $K_{\tS}$ 
is {\em very} ample, and explaining how to show this.
The argument below directly generalizes Catanese's argument 
\cite[p.~502]{catanese} for bidouble covers.

{\noindent {\em Proof of Theorem~\ref{mufti}.  }}
By Proposition~\ref{pardinirocks}(iv), as $\pi$ is finite,\lremind{mulder}
\begin{equation}
H^0 \left( \tS, \cK_{\tS} \right) \cong \oplus_\chi H^0 \left( S, \cK_S \left( L_\chi \right) \right)
\label{mulder}
\end{equation}
We will need to understand this isomorphism more precisely, in
particular how summands on the right of \eqref{mulder} give global
differentials on $\tS$.  For example, the map $H^0(S, \cK_S)
\hookrightarrow H^0(\tS, \cK_{\tS})$ (corresponding to the summand
$\chi=0$) is the pullback map; the pullback of a nonzero $s \in H^0(S, K_S)$
vanishes on the
pullback of the divisor of zeros of $s$, along with the ramification
divisor $R = \oplus_\si R_\si$ with multiplicity $p-1$.  Let $z_\si
\in H^0(\tS, R_\si)$ ($\si \in G$) be a section with divisor $R_\si$
(preserved by the Galois group $G$).
A local calculation gives \lremind{muppet}
\begin{equation} \label{muppet}
H^0 \left( \tS, K_{\tS} \right) \cong \bigoplus_\chi \left( \prod_\si
  z_{\si}^{p-1-\langle \si, \chi \rangle} \right) H^0 \left( S, \cK_S \left( L_\chi \right) \right).
\end{equation}
(By the hypotheses of
Proposition~\ref{pardinirocks}, no more than two $R_\si$ pass through any
point.  We consider three cases.
Case 0: This is clear away from points of $R = \oplus R_\si$.  Case
1: To do this local calculation near points lying on precisely one
$R_{\si}$, use the fact that there are local coordinates $(x,y)$
such that $\pi$ is given by $(x,y)
\rightarrow (x^p, y)$.  Case 2: Near points lying on precisely two
$R_{\si}$, the morphism is given by $(x,y)
\rightarrow (x^p, y^p)$ for appropriate $x$ and $y$.)

We first show that the canonical system $\left| \cK_{\tS} \right|$ is
base point free.  Given a point $q \in \tS$, $\pi(q)$ lies on at most
two $D_\si$.  Choose a $\chi$ such that $\langle \si, \chi \rangle
=p-1$ for all such $\si$; such a $\chi$ exists as $G$ has dimension
$3$ over $\F_p$.  Choose a section of $\cK_S(L_\chi)$ not vanishing at
$\pi(q)$ (possible by sufficient ampleness of $L_\chi$).  Then by
\eqref{muppet} the corresponding section of $\cK_{\tS}$ does not
vanish at $q$.

We next show that $\left| \cK_{\tS} \right|$ separates
points.  Because $\cK_S(L_\chi)$ separates points for any $\chi \neq 0$ 
(by sufficient ampleness of $L_\chi$), $\left| K_{\tS} \right|$ separates
points separated by $\pi$.   Suppose now that
$\pi(p_1)=\pi(p_2)$. 
For each $\chi \neq 0$, choose a section $s_\chi$ of $\cK_{S}(L_\chi)$ not
vanishing at $\pi(p_1)$.  The corresponding $\left| G \right|-1$
sections of $\cK_{\tS}$
give a map near $p_1$ and $p_2$ to 
$\proj^{\left| G \right| -2}$ (that
factors through $\left| \cK_{\tS} \right|$).  As described above, an element $g$ of the
Galois group $G$ acts on the section corresponding to $\chi$ 
by the character $\chi$, i.e.\ by
multiplication by the root of unity $\zeta^{\langle g, \chi \rangle}$.
Suppose
that $g(p_1)=p_2$.  If $p_1$ and $p_2$ are mapped to the same point of
$\proj^{\left| G \right| -2}$, then 
$$
\left(  \left( \prod_\si
  z_{\si}^{p-1-\langle \si, \chi \rangle} \right) s_\chi \right)_{\chi \neq 0}
\text{ and }
\left(
\zeta^{\langle g, \chi \rangle}
   \left( \prod_\si
  z_{\si}^{p-1-\langle \si, \chi \rangle} \right) s_\chi \right)_{\chi \neq 0}
$$
are linearly dependent, so 
$\langle g, \chi_1 \rangle =
\langle g, \chi_2 \rangle$ for all $\chi_i$ satisfying $\prod_\si
z_{\si}^{p-1-\langle \si, \chi_i \rangle} \neq 0$ and $\chi_i \neq 0$.
We again have three cases.
Case 0: If no $z_\si$ is zero (i.e.\ the $p_i$ lie in the \'etale locus
of $\pi$), this forces $g$ to be the identity, so $p_1=p_2$.  Case 1: if
exactly one $z_{\si}$ is $0$,  then $\si$ preserves $p_1$.  Now
$\langle g, \chi_1 \rangle = \langle g, \chi_2 \rangle$ for all
$\chi_i$ with $\langle \si, \chi_i \rangle = p-1$.  This set is
non-empty; by translating all such $\chi_i$ by one fixed such $\chi$,
we have that $\langle g, \chi_1 \rangle = \langle g, \chi_2 \rangle$
for all $\chi_i$ with $\langle \si, \chi_i \rangle =0$, i.e.\
$\langle g, \chi \rangle = 0$
for all $\chi$ with $\langle \si, \chi \rangle =0$. By linear
algebra over $\F_p$, $g$ must lie in the subspace generated by
$\si$, i.e.\ $g$ is a multiple of $\si$.  Thus $p_2 = g(p_1)=p_1$.  Case
2: suppose two $z_{\si}$ vanish at $p_1$, say $z_{\si_1}$, $z_{\si_2}$.  Then
for all $\chi_1$, $\chi_2$ in the non-empty set $\{ \chi : \langle
\si_1, \chi \rangle = \langle \si_2, \chi \rangle = p-1 \}$, $\langle
g, \chi_1 \rangle = \langle g, \chi_2 \rangle$.  By translating
by one such $\chi$ we have that
for all $\chi$ such that   $\langle
\si_1, \chi \rangle = \langle \si_2, \chi \rangle = 0$, $\langle
g, \chi \rangle = 0$.  
Again, by linear
algebra on $(\F_p)^3$, $g$ must lie in the subspace generated by
$\si_1$ and $\si_2$, so again 
$p_2 = g(p_1)=p_1$.  

Thus the canonical system $\left| \cK_{\tS} \right|$ separates points.
We conclude the proof by showing that it separates tangent vectors.
Case 0: Near any point disjoint from the $R_i$ (i.e.\ where $\pi$ is
\'etale), the sections of $\cK_{\tS}$ corresponding to $\cK_S(L_\chi)$
separate tangent vectors for any $\chi \neq 0$.  Case 1: Suppose $q$
lies in precisely one $R_{\si}$.  It suffices to exhibit two sections
$s_1$, $s_2$ of $\cK_{\tS}$ vanishing at $q$ to precisely first order,
whose tangent directions are transverse.  First choose $\chi$ such
that $\langle \si, \chi \rangle =p-1$, and a section of
$\cK_S(L_\chi)$ vanishing to order $1$ at $\pi(q)$, whose zero-set is
transverse to $D_\si$ at $\pi(q)$; then the corresponding section
$s_1$ of $\cK_{\tS}$ vanishes to order $1$ at $q$, and its zero-set is
transverse to $R_\si$.  Second, choose $\chi \neq 0$ such that
$\langle \si, \chi \rangle = p-2$, and a section of $\cK_S(L_\chi)$
not vanishing at $q$; the corresponding section $s_2$ of $\cK_{\tS}$
vanishes to order $1$ at $q$ and its zero-set is contains $R_\si$.
Case 2:  Suppose $q$ lies on $R_{\si_1} \cap R_{\si_2}$.
For $i=1,2$, choose $\chi_i \neq 0$ such that $\langle \si_i, \chi_i \rangle = p-2$
and $\langle \si_{3-i}, \chi_i \rangle = p-1$ (possible as 
$\si_1$ and $\si_2$ are linearly independent in $G \cong (\F_p)^3$,
by the hypotheses of Proposition~\ref{pardinirocks}). 
Then near $q$, $s_i$ vanishes precisely along $R_{\si_i}$.
 \epf

\tpoint{Theorem} {\em  In Examples~\ref{p2} and~\ref{p3}, if $A$ is  sufficiently ample,
then: \label{regularity} \lremind{regularity}
\begin{enumerate}
\item[(a)]
$\tS$ is regular:  $q(\tS) := h^1(\tS, \oh_{\tS}) = 0$.
\item[(b)]
The deformations of $\tS$ are the same as the deformations
of $(S, \{ D_\si \})$.  In particular, the deformations
of $G$-covers are also $G$-covers.
\item[(c)] 
The deformation space  of $\tS$ has the same type as the deformation space 
of $(S, C)$. 
\item[(d)] $\Aut(\tS) \equiv G$ (the only automorphisms of $\tS$
are those preserving the cover of $S$).
\end{enumerate}
}

Part (c) implies that the fine moduli stack
of surfaces of general type satisfies Murphy's law.
Part (d) implies that the fine moduli stack
is locally a quotient of the moduli space of $(S, \{ D_\si \}
)$ by a trivial $G$-action (automorphism
groups are semicontinuous in families, see for example
\cite[Corollary~4.5]{fp1}), so 
the coarse moduli space also satisfies Murphy's
Law.  Thus the Proposition implies {\bf \M 2a--c}.

\bpf 
(a) By the Leray spectral sequence, 
$$h^1(\tS, \oh_{\tS}) = 
h^1(S, \pi_* \oh_{\tS}) =
\sum_{\chi} 
h^1 \left( S,  L_{\chi}^{-1} \right) = 0$$
using Serre vanishing (for $\chi \neq 0$) and  
the regularity of any blow-up of $\proj^2$ (for $\chi=0$).

(b)  For example~\ref{p2} ($p=2$),
the result follows from \cite[Corollary~3.23]{manetti};
we restate the three hypotheses of Manetti's result 
for the reader's convenience.
(i)  $S$ is smooth of dimension $\geq 2$, and $H^0(S, T_S)=0$.
(The latter is true because $S$ has no non-trivial infinitesimal automorphism.
Reason:  any such would descend to an infinitesimal automorphism
of $\proj^2$ fixing the $p_i$, in particular $p_1 = [1;0;0]$,
\dots, $p_4 = [1;1;1]$.)
(ii) $H^0(S, T_S(-L_\chi))=\Ext^1_{\oh_S}(\Omega^1_S, L^{-1}_\chi)=
H^1(S, L^{-1}_\chi)=0$ (true by Serre vanishing, and sufficient
ampleness of $A$).
(iii) $H^0(S, D_{\si} - L_\chi)=0$ for all $\chi \neq 0$,
$\langle \si, \chi \rangle = 0$
(true by Serre vanishing).  Hence (b) holds for Example~\ref{p2}.

The paper \cite{manetti} deals with $(\zed/2)^r$ covers.
However, \cite[Corollary~3.23]{manetti} applies without
change for $(\zed/p)^r$-covers.
The only change in the proof arises in the proof of the prior
result \cite[Proposition~3.16]{manetti}; the statement of this
proposition remains the same, and the proof is changed in the obvious way.
In particular, the fourth equation display should read
$$
\Omega^1_{X/Y} = \bigoplus_\si \frac {\oh_X(- (p-1) R_\si)} {w_\si \oh_X ( - p R_\si)}
= \bigoplus_\si \oh_{R_\si}(-R_\si).
$$
Then the hypotheses of \cite[Corollary~3.23]{manetti} follow as in the
case $p=2$, and we have proved (b) for Example~\ref{p3} ($p=3$) as well.

(c) Choose $A = C + n p K_{\tS}$ for $n \gg 0$, so that its higher
cohomology vanishes.  Then $\Def(S,\{ D_\si \}) \rightarrow \Def (S,
C)$ is a smooth morphism: in any deformation of $S$ the divisor class
$[D_{\si}]$ extends (as $C$ and $K_{\tS}$ extend), and extends uniquely
(by $h^1(S, \oh_S)= 0$), and the choice of divisor in the divisor
class is a smooth choice. 

(d) follows from \cite[Theorem~4.6]{fp1} (taking $D_1$ of \cite{fp1}
to be any of the $D_{\si}$ in class $A$; Fantechi and Pardini's $m_1$ is our $p$).
\epf

\epoint{Proof of {\bf \M 2d}} \label{pf2d} \lremind{pf2d}By taking six
general sections of a sufficiently positive multiple of the canonical
bundle (very ample with vanishing higher cohomology), and using this
to embed $\tS$ in $\proj^5$, we see that the Hilbert scheme of
nonsingular surfaces in $\proj^5$ satisfies Murphy's Law.  (The choice
of the sections, up to scalar, is a smooth choice over the fine moduli
space of surfaces; we use here $h^1(\tS, \oh_{\tS})=0$, by
Theorem~\ref{regularity}(a), so the line bundle over $\tS$ extends
uniquely over the universal surface over the moduli space.)

Using five general sections of the bundle to map $\tS$ to $\proj^4$
yields a surface with only singularities in codimension $2$; each
consists of two nonsingular branches meeting transversely (a
non-Cohen-Macaulay singularity).  The deformations of such a
singularity preserve the singularity.  (This can be checked formally
locally; the calculation can then be done tediously by hand using two transverse
co-ordinate planes in $\A^4$.)  Hence deformations of the singular
surface in $\proj^4$ correspond to deformations of the nonsingular
surface $\tS$ along with the map to $\proj^4$.  We have thus proved
{\bf \M 2d}.

\section{Deformations of products:  proof of {\bf \M 3}} 

We use $\tS$ to construct an example in any dimension $d>2$ of a
$d$-fold with very ample canonical bundle with deformation space of
the same type as that of $\tS$, and with automorphism group $G$; {\bf
  \M 3} then follows.  As $q_{\tS} = h^1(\tS, \oh_{\tS}) = 0$, there
are no nonconstant morphisms $\tS \rightarrow C$ to a nonsingular
curve of positive genus. Let $C_1$, \dots, $C_{d-2}$ be general curves
of genus $3$, and let $Y=\prod C_i$.  Now $\Def \prod C_i = \prod \Def
C_i$ (true on the level of first-order deformations by
\cite[Theorem~2.2]{vo}, cf.\ \cite[Exercise~3.33]{hm}; and $\Def C_i$
is unobstructed), so $Y$ is unobstructed. 

We will need the following improvement of a result of Ran.  (In the
smooth holomorphic case, this result is very similar to Horikawa's
\cite[Theorem~8.2]{horikawa}.)  The result is likely known to experts.

\tpoint{Theorem} {\em Let $f: X \rightarrow Y$ be a morphism with $f_*
  \oh_X = \oh_Y$ and $R^1 f_* \oh_X = 0$.  Then $\Def(f: X \rightarrow
  Y) \rightarrow \Def X$ is an isomorphism.}
\label{wicked} \lremind{wicked}

(See \cite[Definition~1.1]{ran} for a definition of $\Def(f:X
\rightarrow Y)$.)  Ran's theorem \cite[Theorem~3.3]{ran} is identical
to this, except with the additional hypothesis that $R^2 f_* \oh_X =
0$, and the weaker conclusion that $\Def(f: X \rightarrow Y)
\rightarrow \Def X$ is smooth.  This proof is simply  a refinement
of his.

\bpf Consider the spectral sequence $\Ext^p_Y(\Omega_Y, R^q f_* \oh_X)
\Rightarrow \Ext^i_f(\Omega_Y, \oh_X)$ \cite[(6)]{ran}; the $E^2$ term
is
$$
\xymatrix{
\Hom_Y (\Omega_Y, R^2 f_* \oh_X) & \Ext_Y^1(\Omega_Y, R^2 f_* \oh_X) 
  & \Ext_Y^2(\Omega_Y, R^2 f_* \oh_X)  \\ 
0 \ar[rrd] & 0 & 0 \\
\Hom_Y(\Omega_Y,  \oh_Y) &
\Ext^1_Y (\Omega_Y,  \oh_Y) &
\Ext^2_Y (\Omega_Y,  \oh_Y) &
}
$$
and hereafter in the spectral sequence the entry $\Ext^i_Y (\Omega_Y, \oh_Y)$
  will not change for $i=0,1,2$.  Hence we conclude that 
\lremind{ranfun}
\begin{equation}\label{ranfun}
\xymatrix{ T^i_Y := \Ext^i_Y (\Omega_Y, \oh_Y) \ar[r] &
 \Ext^i_f(\Omega_Y, \oh_X)}
\end{equation}
 is an isomorphism for $i=0$ and $1$ and an injection for $i=2$. 

By Ran's exact sequence \cite[(2.2)]{ran}
$$
\xymatrix{
& T^0_X \oplus T^0_Y \ar[r] &
  \Ext^0_f(\Omega_Y, \oh_X)\ar[r] &  \\
 T^1_f \ar[r] & T^1_X \oplus T^1_Y \ar[r] &
  \Ext^1_f(\Omega_Y, \oh_X)\ar[r] &  \\
T^2_f \ar[r] & T^2_X \oplus T^2_Y
  \ar[r] & \Ext^2_f(\Omega_Y, \oh_X) ,}$$
we have that $T^1_f
\rightarrow T^1_X$ is an isomorphism and $T^2_f \rightarrow T^2_X$ is
injective, which gives the desired result. 
\epf

\tpoint{Proposition} \label{mummy} \lremind{mummy}
{\em  \begin{enumerate}
\item[(a)] $\Def \left( \tS \times Y \right) \cong \Def \tS \times \Def Y$
\item[(b)]  $\Aut \left( \tS \times Y \right) \cong G$ \end{enumerate}
} 

Part (a) shows that the deformation space of $\tS \times Y$ has the
same type as that of $\tS$, and hence the fine moduli stack satisfies
Murphy's Law.  Part (b) shows that the $d$-fold has no ``extra''
automorphisms, and thus (as in the surface case) the coarse moduli space
satisfies Murphy's Law.  Hence {\bf \M 3} follows.

\bpf (a) We first show that the natural morphism $\Def \tS \times \Def
Y \rightarrow \Def (\tS \times Y \rightarrow Y)$ is an isomorphism.
For convenience let $X = \tS \times Y$.  
Let 
\begin{equation} 
\label{mugwump}
\xymatrix{
\cX \ar[rr] \ar[dr] & &  \cY  \ar[dl] \\
& \Def(X \rightarrow Y)
}
\end{equation}
be the universal morphism over $\Def (X \rightarrow Y)$; all morphisms
in \eqref{mugwump} are flat.  By the flatness of the horizontal
morphism $\cX \rightarrow  \cY$ of \eqref{mugwump} we have a
natural morphism $\cY \rightarrow \Def \tS$.  This morphism
descends to $\Def(X \rightarrow Y) \rightarrow \Def \tS$.  (Reason:
Interpret 
\begin{equation} 
\label{mug}
\xymatrix{
\cY \ar[r] \ar[d] & \Def \tS \\
\Def(X \rightarrow Y) 
}
\end{equation}
as a family of maps from irreducible projective varieties
to $\Def \tS$ over a formal local scheme
$\Def(X \rightarrow Y)$, that is constant on the
central fiber.  As constant maps deform to constant maps, 
\eqref{mug} must factor through $\Def(X \rightarrow Y)
\rightarrow \Def \tS$.)

Hence we have a natural morphism $\Def(\tS \times Y \rightarrow Y) \rightarrow
\Def \tS \times \Def Y$; there is of course a natural morphism
in the other direction.  By observing how the universal families
behave under these morphisms, we see that the morphisms are
isomorphisms (and mutual inverses).

Finally, $\Def \left( \tS \times Y \rightarrow Y \right) \rightarrow
\Def \left( \tS \times Y \right)$ is an isomorphism by
Theorem~\ref{wicked}.  The last hypothesis of Theorem~\ref{wicked}
follows from Proposition~\ref{regularity}(a).
\remind{Old note:  before I improved Ran's result, I had a clunky proof
of this paragraph.
Look for it up to Mnov2004.tex.}

(b) Suppose $X$ and $Y$ are varieties where (i) $Y$ is connected and has
no nontrivial automorphisms, (ii) $X$ has discrete automorphism group, and
(iii) the only morphisms from $X$ to $Y$ are constant morphisms (i.e.\ 
$\Hom(X,Y) = Y$).  Then $\Aut X \rightarrow \Aut(X \times Y)$
is an isomorphism.  ({\em Proof:} Define the {\em fibers} of $X \times
Y$ to be the fibers of the projection to $Y$, so they are
canonically isomorphic to $X$.  Consider any automorphism $a: X \times
Y \rightarrow X \times Y$. Each fiber of the source must map to a
fiber of the target by (iii).  Hence $a$ induces an
automorphism of $Y$, necessarily the identity by (i).
Compose $a$ with an automorphism of $X$ so that $a$ is the identity on
one fiber.  Then $a$ is the identity on all fibers by (ii)
and the connectedness of $Y$.)

As the $C_i$ are chosen generally, the only morphisms $C_i \rightarrow
C_j$ ($i \neq j$) are constant maps.  Recall that 
there are no nonconstant maps $\tS \rightarrow C_i$, and 
$C_i$ has no automorphisms.  
Hence by induction on the factors of $\tS \times (\prod C_i)$,
$\Aut (\tS \times Y) \cong \Aut \tS$. \epf

Pardini points out that in the first paragraph of the proof of (b),
the right way to see that every automorphism of $X \times Y$ sends
a fiber to the fiber is that the projection
$\tS \times C_1 \times \cdots \times C_k \rightarrow 
C_1 \times \cdots \times C_k$ is (essentially) the Albanese map.

\section{Branched covers of $\proj^2$ and the proof of {\bf \M 5}}

Take three sections of a sufficiently positive bundle on $\tS$.
As in the proof of {\bf \M 2d} (Section~\ref{pf2d}), 
the line bundle over $\tS$ 
\'etale-locally extends uniquely over the universal surface over the moduli
space, and the choice of three sections (up to non-zero scalar)
is a smooth one.  Hence we have {\bf \M 5a}.
J. Wahl provides the connection to {\bf \M 5b}:

\tpoint{Theorem (Wahl \cite[p.~530]{wahl})} 
{\em Let $Y \rightarrow \proj^2$ be a finite surjective
morphism, $Y$ a nonsingular surface, whose branch curve $C$ is
reduced with only nodes and cusps as singularities.  Then via taking
branch curves, there is a one-to-one correspondence between infinitesimal
deformations of the morphism $Y \rightarrow \proj^2$ and infinitesimal
deformations of $C$ in $\proj^2$ which preserve the formal
nature of the singularities.} \lremind{wahlthm} \label{wahlthm}

Wahl's paper assumes that the characteristic is $0$, but his proof of
this result uses only that the characteristic is not $2$ or $3$.  
To reassure the reader, we point out the places where
characteristic $0$ is used before Wahl's proof of
Theorem~\ref{wahlthm}
concludes on p.~558.  Proposition 1.3.1 and equation (1.5.3)
are not used in the proof.  Theorem 2.2.8 and its rephrasing
(Theorem 2.2.11) give
a normal form for stable singularities, and use only that the
characteristic is not $2$ or $3$.
(One might conjecture that an appropriate  formulation is true in characteristic $2$
and $3$, but I have not attempted to prove this.)
Part {\bf \M 5b} then follows from the next result.

\tpoint{Proposition} {\em If $\tS$ is any smooth projective surface over an
  infinite base field of characteristic not $2$ or $3$, and $\cL'$ is
  an ample invertible sheaf, then for $n \gg 0$, three general
  sections of ${\cL'}^{\otimes n}$ give a morphism to $\proj^2$ with
  reduced branch curve with only nodes and cusps as singularities.}
\lremind{nicebranch} \label{nicebranch}

In characteristic $0$ the result is classical (presumably nineteenth
century); the proof is by taking $n$ large enough that
${\cL'}^{\otimes n}$ is very ample, and then taking a generic
projection.  Because we need the result in positive characteristic as
well, we use a slightly different approach, although as usual we show
the result by showing that ``nothing worse can happen,'' by excluding
possibilities on a case-by-case basis.

\bpf We will make repeated use of the following useful fact
\cite[Example~12.1.11]{it} without comment: Let $E$ be a vector bundle on
a variety $X$ over an infinite field, generated by a
finite-dimensional vector space $W$ of sections.  Let $V$ be a
subvariety of (the total space of) $E$ of pure codimension $m$.  Then for a general element
of $W$, the pullback of $V$ to $X$ has pure codimension $m$.

We show that the branch curve is as desired for maps to $\proj^2$
given by  three general sections
of $\cL$, for all invertible sheaves $\cL$:
\begin{itemize}
\item  separating $3$-jets (i.e.\ 
$J^3 (\cL)$ is generated by global sections of $\cL$), 
\item separating
$2$-jets at pairs of distinct points (i.e.\ $\pi_1^* J^2 (\cL) \oplus
\pi_2^* J^2 (\cL)$ on $\tS \times \tS - \Delta$ is generated by
global sections of $\cL$ on $\tS$), 
\item and separating $1$-jets at
triples of distinct points.
\end{itemize}
A sufficiently high power of
$\cL'$ certainly has these properties. 

We first show that {\em the ramification locus is codimension $1$.}
The rank $9$ bundle ($3 \times 3$ matrix bundle)
$$E_1=\Hom(\oh_{\tS}^{\oplus 3}, J^1(\cL)) \cong J^1(\cL)^{\oplus 3}$$
is generated by global sections of $\cL$.  The rank $\leq 2$
(determinant $0$) locus of $E_1$ is codimension $1$.  Thus the
ramification locus on $\tS$, i.e.\ where the induced section of $E_1$ is
rank $\leq 2$, indeed has pure codimension $1$.  The locus where the
section of $E_1$ meets the rank $\leq 1$ locus is codimension $4$ on
$\tS$, i.e.\ empty.  Hence at each point $p$ of the ramification curve,
the section of $E_1$ has rank exactly $2$, so there is a section $x$
in our net vanishing at the point, but not vanishing to second order
($x$ is a local coordinate on $\tS$).  There is another section $z$
corresponding to the kernel $\cJ := \ker(\oh_{\tS}^{\oplus 3} \rightarrow
J^1(\cL))$ that vanishes to order (at least) $2$ at $p$.  (Note that $x$ and
$z$ are lines on $\proj^2$ and hence local coordinates there.)  The
section $z$ is unique up to multiplication by nonzero scalar. The
section $x$ is unique up to multiplication by nonzero scalar, and
addition by some multiple of $z$. \remind{Old note:  Reverse description of $x$
and $z$?}

Consider next the $6 \times 3$ matrix bundle
$E_2=\Hom(\oh_{\tS}^{\oplus 3}, J^2(\cL))$.  
The ``rank $2$ locus of $E_1$'' makes sense for $E_2$ (and indeed 
any jet bundle 
surjecting onto $E_1$).
On this subvariety of $E_2$ there is a map of invertible sheaves $\al:
\cJ \rightarrow 
(\Omega^2_{\tS} \otimes \cL) / (x \Omega^1_{\tS})$ locally 
induced by\lremind{mutton}
\begin{equation}
\label{mutton}
\xymatrix{ 0 \ar[r] & \cJ \ar[r] \ar[d]^{\al} & \oh_{\tS}^{\oplus 3} \ar[d]
  \\
  0 \ar[r] & \frac{ \Omega^2_{\tS} \otimes \cL}  {x \Omega^1_{\tS}} \ar[r]
  & \frac{ J^2 (\cL)} { x J^1(\oh_{\tS})} \ar[r] & \frac{J^1(\cL) }{ x \oh_{\tS} }\ar[r] & 0.  }\end{equation}
Note that both
$(\Omega^2_{\tS} \otimes \cL) / (x \Omega^1_{\tS})$ and 
$\al$ are independent of $x$.
(Near a point $p$ of the ramification curve the bottom exact sequence
has a straightforward interpretation.  If $y$ is a local coordinate of
$\tS$ at $p$ transverse to $x$, then once a trivialization of $\cL$ near $p$
is chosen, the left term corresponds to the coefficient of $y^2$, and
the right term corresponds to coefficients of $1$ and $y$.)  Let $V$
be the subvariety of $E_2$ 
where $\al$  has rank $0$
(inside the rank $2$ locus of $E_1$); then $\codim V=2$, so the corresponding
locus on $\tS$ has codimension $2$ as well.  For the purposes of this
proof, call such points of the ramification curve {\em twisty} points.
Thus for a non-twisty ramification point $p$, there is a section not
vanishing at $p$, another section $x$ vanishing to first order at $p$,
and a third section $z$ (corresponding to $\cJ$) vanishing to order
{\em exactly} $2$ at $p$ such that $(x,z)$ has length $2$ on $\tS$: 
we conclude that {\em the
  ramification is simple}, and the ramification curve is smooth at $p$.
(Note that $x$ and $z$ are local coordinates on $\proj^2$, and there
are local coordinates $x$ and $y$ on $\tS$ such that formally locally
$z-y^2 \in (y^3, xy, x^2)$.)
We remark that we have defined a codimension $2$ subvariety of $E_2$
(the ``twisty subvariety''), and of any jet bundle surjecting onto $E_2$.

Consider now the locus in the $9 \times 3$ matrix bundle
$$E_3=  \Hom \left(  \oh_{\tS}^{\oplus 3}, 
\pi_1^* J^1(\cL) \oplus \pi_2^*
J^1(\cL) \oplus \pi_3^*
J^1(\cL) \right)$$
on $\tS \times \tS \times \tS - \Delta$
where the image in the $3 \times 3$ matrix bundle
$$
\Hom \left( \oh_{\tS}^{\oplus 3}, \pi_1^* \cL \oplus \pi_2^* \cL \oplus
\pi_3^* \cL \right)$$
has rank $1$ (``the three points of $\tS$ map to the
same point of $\proj^2$,'' codimension $4$ in $E_3$) and where the
image in each of the $3 \times 3$ matrix bundles $\pi_i^* E_1$ has
rank $2$ (``each point is on the ramification curve,'' codimension $1$
in $E_3$).  The pullback of this locus (via our three sections of $\cL$) 
has codimension $4+1+1+1=7$ on $\tS \times \tS
\times \tS - \Delta$, and hence is empty.  Hence {\em there cannot be
  three points of the ramification curve mapping to the same point of
  $\proj^2$.}

Consider the locus in  the $9 \times 3$ matrix bundle
$$E_4= \Hom \left( \oh_{\tS}^{\oplus 3}, \pi_1^* J^2(\cL) \oplus \pi_2^*
J^1(\cL) \right)$$
on $\tS \times \tS - \Delta$, where we require: the image
in the $2 \times 3$ matrix bundle $\Hom(\oh_{\tS}^{\oplus 3}, \pi_1^*
\cL \oplus \pi_2^* \cL)$ to be rank $1$ (``the two points map to the
same point of $\proj^2$,'' a codimension $2$ condition); the first
point to be a twisty ramification point (shown earlier to be a
codimension $2$ condition on $E_4$); and the second
point to be a ramification point (a codimension $1$ condition on 
$E_4$ as shown above).  These conditions are 
independent, so this locus in $E_4$ has codimension $5$, hence the
pullback (by our section) is empty on $\tS \times \tS - \Delta$.  Thus
{\em if two points of the ramification curve map to the same point of
  $\proj^2$, neither is twisty.}

Next, in the $6 \times 3$ matrix bundle
$$E_5=\Hom \left( \oh_{\tS}^{\oplus 3}, \pi_1^* J^1(\cL) \oplus \pi_2^*
J^1(\cL) \right)$$
on $\tS \times \tS - \Delta$, we consider the locus where: the
image in the $2 \times 3$ matrix bundle $\Hom(\oh_{\tS}^{\oplus 3},
\pi_1^* \cL \oplus \pi_2^* \cL)$ has rank $1$ (``the two points map to
the same point in $\proj^2$,'' codimension $2$); and $E_5$ itself 
has rank $\leq 2$ (``the two points are ramification
points, and the branch curves in $\proj^2$ share a tangent line,'' contributing
an additional codimension of $3$ --- not $4$, because of 
the dependence of these
 conditions with the previous ones).  The total codimension of this
locus is $5$, hence (the pullback of) 
this locus is empty on $\tS \times \tS - \Delta$,
so {\em if two ramification points map to the same point of $\proj^2$,
  the tangent vectors of the branch curve in $\proj^2$ are
  transverse}.  In particular, {\em the branch curve is reduced,
and has only nodes away from the twisty points}.

Finally we  show that the twisty ramification points give
cusps of the branch curve.  Similar to \eqref{mutton}, on the
ramification curve we locally have
a morphism $\be$ from $\cJ=\ker (\oh_{\tS}^{\oplus 3} \rightarrow J^1(\cL))$
 to a rank $2$ bundle $Q$ satisfying defined by
$$
\xymatrix{ 0 \ar[r] & \cJ \ar[r] \ar[d]^{\be} & \oh_{\tS}^{\oplus 3} \ar[d]
  \\
  0 \ar[r] & Q  \ar[r]
  & \frac{ J^3 (\cL)} { x J^2(\oh_{\tS})} \ar[r] & \frac{J^1(\cL) }{ x \oh_{\tS} }\ar[r] & 0.  }
$$
satisfying
$$
\xymatrix{
& & \cJ \ar[d]^\be \ar[dr]^\al \\
0 \ar[r] &  \frac { \Omega^3_{\tS} \otimes \cL} { x \Omega^2_{\tS}} 
\ar[r] & Q \ar[r]  & 
\frac { \Omega^2_{\tS} \otimes \cL} { x \Omega^1_{\tS}} 
\ar[r] & 0.}$$
(Near a point $p$ of the ramification curve, the local
interpretation is similar to that of \eqref{mutton}.  The term
${( \Omega^3_{\tS} \otimes \cL )} / { x \Omega^2_{\tS}}$
corresponds to the coefficient of $y^3$, and the
term ${ (\Omega^2_{\tS} \otimes \cL)} / { x \Omega^1_{\tS}}$
corresponds to the coefficient of $y^2$.)  Hence on the
(codimension $2$) twisty locus, where $\al =0$, there is a morphism
of invertible sheaves $\ga:
\cJ \rightarrow (\Omega^3_{\tS} \otimes \cL ) / (x \Omega^2_{\tS})$; and
both
$(\Omega^3_{\tS} \otimes \cL ) / (x \Omega^2_{\tS})$ and $\ga$ are
independent of our choice of $x$.  Consider the
subvariety of $\Hom(\oh_{\tS}^{\oplus 3}, J^3(\cL))$ (in the
codimension $2$ twisty locus) such that the rank of $\ga$ is $0$.
This is a codimension $3$ subvariety of $\Hom(\oh_{\tS}^{\oplus 3},
J^3(\cL))$, and hence empty on $\tS$.  Thus formally locally, at a
twisty ramification point, ``the $y^3$-coefficient is non-zero,'' 
so the morphism is given by
$$(x,y) \mapsto \left( x, y^3 f_1(x,y) + y^2 (x g_1(x)) + y (x h_1(x))
  + x^2 i_1(x) \right)=(x,z)$$
where $f_1(0,0) \neq 0$.  Consider
again the $6 \times 3$ matrix bundle $E_2=\Hom( \oh_{\tS}^{\oplus 3},
J^2(\cL))$.  The rank $2$ locus of $E_2$ is codimension $4$, and hence
our section misses it: every non-zero section of $\cL$ restricts to
something nonzero in $J^2(\cL)$ at every point. Hence $z \notin 
(x,y)^3$, so either $i_1(0) \neq 0$ or $h_1(0) \neq 0$. By
replacing $x$ by $x+y$ if necessary, we may assume that $h_1(0) \neq
0$.

By replacing $y$ by a scalar multiple, we may assume $f_1(0,0)=1$.
By replacing $y$ by $y f_1(x,y)^{1/3}$,  and rearranging, the morphism
may be rewritten as
$$(x,y) \mapsto \left( x, y^3(1+x f_2(x,y)) + y^2 (x g_2(x)) + y (x h_2(x)) + x^2 i_2(x) \right)=(x,z)$$
where again $h_2(x) \neq 0$.
(The $x f_2(x,y)$ arose because of contributions of $x ( y^2 g_1(x) + y h_1(x))$
in the change of variables.) By replacing $y$ by $y (1+xf_2(x,y))^{1/3}$
we obtain
$$(x,y) \mapsto \left( x, y^3(1+x^2 f_3(x,y)) + y^2 (x g_3(x)) + y (x h_3(x)) + x^2 i_3(x) \right)=(x,z).$$
By repeating this process inductively, and noting that the lower degree
terms of $g_n(x)$, $h_n(x)$, and $i_n(x)$ stabilize, we obtain
$$(x,y) \mapsto \left( x, y^3 + y^2 (x G_1(x)) + y (x H_1(x)) + x^2 I_1(x) \right)=(x,z)$$
where $H_1( 0) \neq 0$.
Replacing $y$ by $(y+xG_1(x))/3$, the morphism may be rewritten as
$$(x,y) \mapsto \left( x, y^3 +  y (x H_2(x)) + x^2 I_2(x) \right)=(x,z).$$
Replacing $z$ by $z-x^2 I_2(x)$ and then replacing $x$
by $x H_2(x)$ (here finally using $H_2(0) \neq 0$),
we 
have shown that near a twisty point in formal local co-ordinates the morphism
is given by 
$(x,y) \mapsto (x,z)=(x, y^3  + xy)$.  
Then the branch locus in the $(x,z)$-plane is given by $4x^3+27z^2=0$, i.e.\ it is
a cusp.
\epf

\section{From surfaces to the rest of Theorem~\ref{mainthm}}

\bpoint{Proof of {\bf \M 1}}
We are fortunate that Fantechi and Pardini have proved precisely the
result that we need for the
proof of {\bf \M 1}.  If $X \subset \proj^n$ is a subscheme,
let $\Hilb(X)$ be the (connected component of the) 
Hilbert scheme containing $[X]$.

\tpoint{Theorem} {\em (a) (Fantechi-Pardini
  \cite[Proposition~4.2]{fp2}) Let $\tS \subset \proj^n$ be a smooth,
  regular, projectively normal surface.  Let $H$ be a smooth
  hypersurface of degree $l$ in $\proj^n$ meeting $\tS$ transversely
  along a curve $C$, and let $U \subset \Hilb(\tS) \times \Hilb(H)$ be
  the open set of pairs $(\tS', H')$ such that $\tS'$ and $H'$ are smooth
  and transverse and $\tS'$ is projectively normal.  If $l \gg 0$, then
  the morphism $U \rightarrow \Hilb(C)$ (induced by the intersection) is
  smooth. \label{murray} 

(b)  Furthermore, $C$ is embedded by a complete linear system.}
\lremind{murray}

We note that Fantechi and Pardini's proof of (a) invokes Kodaira
vanishing to show that if $F$ is a hypersurface of degree $l$ then
$H^1(F, N_{F/ \proj^n})=0$, but this may be easily checked directly, so their
result is not characteristic-dependent.

\noindent {\em Proof of (b).}
If $\cI_{C/\tS}$ is the ideal sheaf of $C$ in $\tS$, we have the exact
sequence 
$$
0 \rightarrow \cI_{C/\tS}(1) \rightarrow \oh_{\tS}(1) \rightarrow 
\oh_C(1) \rightarrow 0.$$
As $\cI_{C/\tS} \cong \oh_{\tS}(-l)$, $h^1(\tS, \cI_{C/\tS}(1))=0$
by Serre vanishing.
Thus $H^0(\tS, \oh_{\tS}(1)) \rightarrow H^0(C, \oh_C(1))$ is surjective.
As $H^0(\proj^n, \oh_{\proj^n}(1)) \rightarrow H^0(\tS, \oh_{\tS}(1))$ is also surjective
($\tS$ is embedded by a complete linear system), the result follows.
\epf

We now prove {\bf \M 1}.
 Choose a sufficiently ample line bundle
on $\tS$, so that the corresponding embedding $\tS \hookrightarrow
\proj^n$ (by the complete linear system) is projectively normal, and
so that the line bundle has no higher cohomology.  The deformation
space of $\tS \hookrightarrow \proj^n$ is smooth over the deformation
space of $\tS$, as described in Section~\ref{pf2d}.  Then Theorem~\ref{murray}(a) gives {\bf \M 1a}.
Deformations of a smooth curve in $\proj^n$ are the same as
deformations of the corresponding immersion, yielding {\bf \M 1b}.
Theorem~\ref{murray}(b) gives {\bf \M 1c}.

\bpoint{Proof of {\bf \M 4}} Near a seminormal point of the Hilbert
scheme, there is a morphism from the Hilbert scheme to the Chow variety
\cite[Theorem~6.3]{kollar}.  If the point of the Hilbert scheme
parametrizes an object that is geometrically reduced, normal, and of
pure dimension, then this morphism is a local isomorphism
\cite[Corollary~6.6.1]{kollar}.  Hence {\bf \M 4} follows from {\bf \M 1a}
and {\bf \M 2d}.

\bpoint{Proof of {\bf \M 6}} (I am grateful to R.~Thomas for
explaining how to think about this problem, and for greatly shortening
the following argument.)  The sheaf in question will be the ideal
sheaf $\cI$ of $\tS$ embedded in $\proj^5$ (from {\bf \M 2d}).
The next result 
implies {\bf \M 6}. 

\tpoint{Proposition} {\em If $Y$ is a nonsingular variety with $h^1(Y,
  \oh_Y) = h^2(Y, \oh_Y)=0$, and $X \hookrightarrow Y$ is a local
  complete intersection, then the deformation space of $X \hookrightarrow Y$ is
 canonically isomorphic to  the deformation space of the ideal sheaf $\cI$ of $X$.}

{\em Warning:} This result does not hold for general $X \hookrightarrow Y$.

\bpf
As usual, we describe an isomorphism of first-order deformations
and an injection of obstructions.  
The first-order deformations and  obstructions
for $X \hookrightarrow Y$ are
$H^0(X, N_{X/Y})$ and $H^1(X,N_{X/Y})$ respectively.  The
first-order deformations and obstructions of the torsion-free
sheaf $\cI$ are
$\Ext^1(\cI, \cI)$ and $\Ext^2(\cI, \cI)$ respectively.
The $E^2$ term 
of the local-to-global spectral sequence for 
$\Ext^{\cdot}(\cI, \cI)$
is:
$$
\xymatrix@=12pt{
\vdots & \vdots & \vdots & \\
H^0(Y,\sExt^1(\cI,\cI))\ar[drr] & H^1(Y, \sExt^1(\cI,\cI)) & H^2(Y, \sExt^1(\cI,\cI)) & \cdots \\
H^0(Y,\sHom(\cI,\cI)) & H^1(Y, \sHom(\cI,\cI)) & H^2(Y, \sHom(\cI,\cI)) & \cdots. }
$$
A straightforward and well-known argument yields $\sExt^q(\cI,
\cI) \cong \wedge^q N_{X/Y}$, so $E^2$ may be written as
$$
\xymatrix@=12pt{
\vdots & \vdots & \vdots & \\
H^0(X, N_{X/Y}) & H^1(X, N_{X/Y}) & H^2(X, N_{X/Y}) & \cdots \\
k & 0 & 0 & \cdots }
$$
using $h^i(Y, \oh_Y)=0$ for $i=1$, $2$.
Thus we have $H^0(X,N_{X/Y}) \cong \Ext^1(\cI, \cI)$ and $H^1(X, N_{X/Y})
\hookrightarrow \Ext^2(\cI, \cI)$, concluding the proof.
\epf

\bpoint{Proof of {\bf \M 7}} We obtain the threefold singularity by
embedding our surface in projective space by a complete linear system
arising from a sufficiently positive line bundle (as in the proof of
{\bf \M 2d}).  The deformations of the cone over the surface
are the same as the deformations of the surface in projective
space, by the following theorem of Schlessinger.

\tpoint{Theorem (Schlessinger \cite[Theorem~2]{sch})}
{\em Let $\tS \subset \proj^n$ be a projectively normal variety
(over a field) of dimension $\geq 2$, such that 
$$
h^1 \left( \tS, \oh_{\tS}(v) \right)=
h^1 \left( \tS, T_{\tS}(v) \right)= 0$$
for $v>0$.  Then the versal deformation spaces of $\tS$ in
$\proj^n$ and the singularity $C_{\tS}$ (the cone over $\tS$) are isomorphic.
}

(Although Schlessinger works in the complex analytic category, his
proof is purely algebraic, and characteristic-independent.)
This singularity is Cohen-Macaulay by the following result, concluding
the proof of {\bf \M 7}.

\tpoint{Proposition} {\em Suppose $\tS$ is a Cohen-Macaulay
  scheme (over a field), $h^i(\tS,\oh_{\tS})=0$ for $i=1$, \dots, $\dim \tS - 1$ and
  $h^0(\tS, \oh_{\tS}) = 1$.  Then the embedding of $\tS$
  by a sufficiently ample line bundle is arithmetically Cohen-Macaulay.}

This result follows from a statement of Hartshorne and Ogus
\cite[p.~429 \#3]{ho}.  See \cite[p.~207--8]{gk} or
\cite[Lemma~1.1(2)]{ch} for a proof.  The hypotheses follow from the
regularity of $\tS$, Theorem~\ref{regularity}(a).
(It turns out that in characteristic $0$, $2 K_{\tS}$ is ample
enough, using Kodaira vanishing.)
 
} 

\end{document}